\documentclass[twoside,11pt,reqno]{amsart}
\usepackage{amsmath,amssymb,amscd,mathrsfs}

\makeatletter

\@settopoint\oddsidemargin
\@settopoint\marginparwidth
\@settopoint\evensidemargin
\@settopoint\topmargin
\@addtoreset{equation}{section}

\def\Point#1{\addtocounter{section}{1}\vspace{2mm}
\noindent\S\thesection. {\bf #1.}\def\@currentlabel{\thesection}\setcounter{equation}{0}}

\newtheorem{Lemma}[equation]{Lemma}
\newtheorem{Theorem}[equation]{Theorem}

\def\ocirc#1{\stackrel{_{\,\circ}}{#1}}

\def\ad{{\operatorname{ad}}}

\def\id{\operatorname{id}}

\def\C{{\mathbb C}}

\def\Z{{\mathbb Z}}
\def\Mtype{\mathtt{M}}
\def\Qtype{\mathtt{Q}}
\def\0{{\bar 0}}
\def\1{{\bar 1}}

\def\hom{{\operatorname{Hom}}}

\def\End{{\operatorname{End}}}

\def\wt{{\operatorname{wt}}}

\def\underbar{\mathpalette\@underbar}
\def\@underbar#1#2{\settowidth{\@tempdimb}{$#1#2$}\@tempdimb=0.8\@tempdimb
                   \ooalign{$#1#2$\crcr%
                         \hfil\rule[-.5mm]{\@tempdimb}{.4pt}\hfil}}

\def\bi{{\underbar{i}}}
\def\bj{{\underbar{j}}}

\def\eps{{\varepsilon}}
\def\phi{{\varphi}}

\def\la{{\lambda}}

\def\al{{\alpha}}

\newdimen\hoogte    \hoogte=12pt    
\newdimen\breedte   \breedte=14pt  
\newdimen\dikte     \dikte=0.5pt 

\newenvironment{Young}{\begingroup
       \def\vr{\vrule height0.89\hoogte width\dikte depth 0.2\hoogte}
       \def\fbox##1{\vbox{\offinterlineskip
                    \hrule height\dikte
                    \hbox to \breedte{\vr\hfill##1\hfill\vr}
                    \hrule height\dikte}}
       \vbox\bgroup \offinterlineskip \tabskip=-\dikte \lineskip=-\dikte
            \halign\bgroup &\fbox{##\unskip}\unskip  \crcr }
       {\egroup\egroup\endgroup}
\def\diagram#1{\relax\ifmmode\vcenter{\,\begin{Young}#1\end{Young}\,}\else%
              $\vcenter{\,\begin{Young}#1\end{Young}\,}$\fi}

\begin{document}
\title[Cartan determinants]{\boldmath
Cartan determinants and
Shapovalov forms}
\author{\sc Jonathan Brundan and Alexander Kleshchev}
\address
{Department of Mathematics\\ University of Oregon\\
Eugene\\ OR~97403, USA}
\email{brundan@darkwing.uoregon.edu, klesh@math.uoregon.edu}

\thanks{Second author
partially supported by the NSF (grant no. DMS-9900134).}
\begin{abstract}
We compute the determinant of the Gram matrix of the Shapovalov form on weight
spaces of the basic representation of an affine
Kac-Moody algebra of ADE type (possibly twisted).
As a consequence, we obtain explicit formulae for the determinants
of the Cartan matrices of $p$-blocks of the symmetric group
and its double cover, and of the associated Hecke algebras at roots of
unity.
\end{abstract}
\maketitle

\section{Introduction}

Let $\mathfrak g$ be an affine Kac-Moody algebra
of type $X_N^{(r)}$ as in the table:
\vspace{2mm}
$$
\begin{array}{|l|l|l|l|l|l|l|l|l|}
\hline
X_N^{(r)}&A_{\ell}^{(1)}&D_{\ell}^{(1)}&E_{\ell}^{(1)}&
A_{2\ell-1}^{(2)}&A_{2\ell}^{(2)}&D_{\ell+1}^{(2)}&E_6^{(2)}&D_4^{(3)}\\
\hline
\ell&\geq 1&\geq 4&6,7\hbox{ or }8&\geq 3&\geq 1&\geq 2&4&2\\
\hline
k&0&0&0&\ell-1&\ell&1&2&1\\
\hline
\alpha&\ell+1&4&9-\ell&2&1&2&1&1\\
\hline
\beta&1&1&1&\ell&2\ell+1&2&3&2\\
\hline
\end{array}
$$
\vspace{1mm}

\noindent
We are interested here in the basic representation $V = V(\Lambda_0)$
of $\mathfrak g$, see \cite{Kac}.
Let $|0\rangle$ be a vacuum vector and define the
{lattice}
$V_\Z := U_\Z |0 \rangle$
in $V$, where $U_\Z$ is the $\Z$-subalgebra of 
the universal enveloping algebra
of $\mathfrak g$ generated by the divided powers
$$
e_i^n / n!, \quad f_i^n / n!\qquad(i = 0,1,\dots,\ell,\ n \geq 1)
$$
in the Chevalley generators.
Let $(.\,,.)_S$ denote the Shapovalov form, 
the unique Hermitian form on $V$ satisfying
$(|0\rangle, |0\rangle)_S = 1$ and $(e_i v, v')_S = (v, f_i v')_S$
for $i = 0,\dots,\ell$ and all $v,v' \in V$.
Its restriction to $V_\Z$ gives a symmetric bilinear form
$$
(.\,,.)_S:V_\Z \times V_\Z \rightarrow \Z.
$$ 
Our Main Theorem gives an explicit formula for the determinant of the Gram matrix of this form on each weight space of $V_\Z$.

To state the result precisely, 
recall the description of the weights of $V$ \cite[$\S$12.6]{Kac}:
every weight is of the form $w\Lambda_0 -d \delta$
for some $w$ in the Weyl group $W$ associated to $\mathfrak g$
and some integer $d \geq 0$.
Also
let $\mathscr P(d)$ denote the set of all partitions
$\la = (\la_1 \geq \la_2 \geq \dots)$ of $d$.
Given $\la \in \mathscr P(d)$, we can gather together its equal parts
to represent it as $\la  = 
(1^{r_1} 2^{r_2} \dots)$.
Also recall the number $r\in\{1,2,3\}$ which comes from the type $X_N^{(r)}$. 
Then:

\vspace{2mm}
\noindent
{\bf Main Theorem.} {\em The determinant of the Gram matrix of the 
Shapovalov form on the 
$(w\Lambda_0 - d \delta)$-weight space of $V_\Z$ is
$\alpha^{a(d)}\beta^{b(d)}$ where
$a(d) = \sum_{\la \in \mathscr P(d)} a_\la$,
$b(d) = \sum_{\la \in \mathscr P(d)} b_\la$ and
for $\la = (1^{r_1} 2^{r_2} \dots)$, 
\begin{align*}
a_\la &=
\prod_{i\text{ with }r|i} \binom{\ell+r_i-1}{r_i}\cdot
\prod_{i\text{ with }r \nmid i} \binom{k+r_i-1}{r_i}\cdot\sum_{i\text{ with }r|i} \frac{r_i}{\ell},\\
b_\la
&=
\prod_{i\text{ with }r|i} \binom{\ell+r_i-1}{r_i}\cdot
\prod_{i\text{ with }r\nmid i} \binom{k+r_i-1}{r_i}\cdot\sum_{i\text{ with }r\nmid i} \frac{r_i}{k},
\end{align*}
$\ell,k,\alpha,\beta$ being as in the above table.
The generating functions
$a(q) = \sum_{d \geq 0} a(d) q^d$ and
$b(q) = \sum_{d \geq 0} b(d) q^d$ are given by the formulae
\begin{align*}
a(q)&=  T(q^r)P(q)^k P(q^r)^{\ell-k},\\
b(q) &= (T(q)-T(q^r))P(q)^k P(q^r)^{\ell-k}
\end{align*}
where $P(q) = \prod_{i \geq 1} \frac{1}{1-q^i}$ is the generating function
for the number of partitions of $d$ and 
$T(q) = \sum_{i \geq 1} \frac{q^i}{1-q^i}$ is the generating function for
the number of divisors of $d$.
}

\vspace{2mm}

In \cite{CKK}, De Concini, Kac and Kazhdan constructed
the basic representation over $\Z$ (at least in the untwisted cases)
using an integral version of the vertex operator construction of \cite{KF}.
They showed in particular that the basic representation
remains irreducible on reduction modulo $p$ if and only if
$p \nmid \det X_N$, where $X_N$ is the Cartan matrix of the underlying
finite root system; this also follows immediately from our Main Theorem
on noting that $\det X_N = \alpha \beta^{r-1}$.

Our interest in the theorem comes instead from modular representation
theory.
Suppose now that $\mathfrak g$ is of type $A_\ell^{(1)}$ and set $p = 
(\ell+1)$.
Let $FS_n$ denote the group algebra of the symmetric group over a field
$F$ of characteristic $p$ (assuming in this case that $p$ is prime), 
and let $H_n$ denote the 
Iwahori-Hecke algebra associated to $S_n$
over an arbitrary field but at a primitive $p$th root of $1$
(this case making sense for arbitrary $p \geq 2$).
By \cite{Ariki, G}, there is an isomorphism between 
the basic representation $V_\Z$ of $\mathfrak g$ and the 
direct sum $K = \bigoplus_{n \geq 0} K_n$ of the Grothendieck groups $K_n$  
of finitely generated projective $FS_n$- (resp. $H_n$-) modules for all $n$.
Under the isomorphism, the weight spaces of $V_\Z$ are in 1--1 correspondence
with the block components 
of $K$, a weight space of the form $w \Lambda - d \delta$
corresponding to a block of $p$-weight $d$ (see e.g. \cite[$\S$5.3]{Mat} 
for the
definition of the $p$-weight of a block).
Moreover, according to \cite[Theorem 14.2]{G}, the Shapovalov form  
corresponds to the usual Cartan pairing
$([P], [Q]) = \dim \hom(P,Q)$ between projective modules $P,Q$.
Thus the theorem has the following immediate corollary:

\vspace{2mm}
\noindent
{\bf Corollary 1.} {\em 
Let $B$ be a block of $p$-weight $d$ of either the
group algebra
$FS_n$ of the symmetric group 
over a field of prime characteristic $p$, or
the Hecke algebra $H_n$ over an arbitrary field but
at a primitive $p$th root of unity, in which case $p\geq 2$ is an arbitrary 
integer.
Then the determinant of the Cartan matrix of $B$ is
$p^{N(d)}$ where $$
N(d) = 
\sum_{\lambda = (1^{r_1} 2^{r_2} \dots) \in \mathscr P(d)} 
\frac{r_1+r_2+\dots}{p-1}
\binom{p-2+r_1}{r_1}
\binom{p-2+r_2}{r_2}
\dots.
$$
The generating function $N(q) = \sum_{d \geq 0} N(d) q^d$
equals $T(q) P(q)^{p-1}$.
}
\vspace{2mm}

It is a classical result
of Brauer that the determinant of the Cartan matrix
of a block of $FS_n$ is a power of $p$ (see \cite[84.17]{CR}).
Donkin \cite{Donkin} has proved similarly
that the determinant of the Cartan matrix of a block of $H_n$
{\em divides} a power of $p$. The corollary shows in particular
that the determinant is exactly a power of $p$, even in those
cases where $p$ is not prime, as had been conjectured by Mathas.
We remark that in the case of blocks of $FS_n$, but not of $H_n$,
the explicit generating function given in the
corollary has also recently been obtained by Bessenrodt and Olsson \cite{BO}
using methods from block theory.

Finally suppose that $\mathfrak g$ is of type $A_{2\ell}^{(2)}$
and set $p = (2\ell+1)$.
In this case, the Main Theorem can be reinterpreted as a
computation of Cartan determinants of the $p$-blocks of
the double covers $\widehat S_n$ of the symmetric group.
Following \cite[$\S$9-c]{BKcrys} for notation,
let $S(n)$ be the twisted group algebra of $S_n$ over an algebraically 
closed field 
$F$ of characteristic $p$ (assuming $p$ is an odd prime in this case), 
and let $W(n)$ be the Hecke-Clifford
superalgebra over an algebraically closed field of
characteristic different from $2$ 
at a primitive $p$th root of unity (for arbitrary odd $p \geq 3$). 
By \cite[7.16, 8.13, 9.9]{BKcrys},
there is an isomorphism between the basic representation $V_\Z$
and the direct sum $K = \bigoplus_{n \geq 0} K_n$ of the Grothendieck
groups of finitely generated projective $S(n)$- (resp. $W(n)$-)
supermodules, under which a weight space of the form
$w \Lambda - d \delta$ maps to a superblock of $p$-bar weight $d$ (see
\cite[$\S$9-a]{BKcrys} for the definition of $p$-bar weight of a superblock),
and the Shapovalov form corresponds to the Cartan pairing on 
projective supermodules (see \cite[$\S$7-c]{BKcrys}). So:

\vspace{2mm}
\noindent
{\bf Corollary 2.} {\em 
Let $B$ be a superblock of $p$-bar weight $d$ of either 
$S(n)$ in odd characteristic $p$, or $W(n)$ at a primitive $p$th root of unity,
in which case $p \geq 3$ is an arbitrary odd integer.
Then the determinant of the Cartan matrix of $B$ is
$p^{N(d)}$ where $$
N(d) = 
\sum_{\lambda = (1^{r_1} 2^{r_2} \dots) \in \mathscr P(d)} 
\frac{2r_1+2r_3+2r_5+\dots}{p-1}
\binom{\frac{p-3}{2}+r_1}{r_1}
\binom{\frac{p-3}{2}+r_2}{r_2}
\binom{\frac{p-3}{2}+r_3}{r_3}\dots.
$$
The generating function $N(q) = \sum_{d \geq 0} N(d) q^d$ equals
$(T(q)-T(q^2))P(q)^{(p-1)/2}$.
}
\vspace{2mm}

It is more natural from the point of view of finite group theory
to ask for the Cartan determinant
of a block $B$ of the twisted group algebra 
$S(n)$ in the usual ungraded sense.
According to Humphreys' classification \cite{Hum}, see also 
\cite[9.16]{BKcrys}, we can associated to $B$ its
$p$-bar weight $d$ and a type $\eps \in \{\Mtype, \Qtype\}$.
In case $\eps = \Mtype$, $B$ coincides with a superblock of $p$-bar weight $d$
and it is immediate that its Cartan determinant is as in Corollary 2.
But in the cases when $\eps = \Qtype$ and $d > 0$, the Cartan matrix
of $B$ has twice as many rows and columns 
as the Cartan matrix of the corresponding superblock. 
Nevertheless, we believe the Cartan determinant is the same, based on
explicit computations for small $d$.
In other words, we conjecture that
Cartan determinants of $p$-blocks of $S(n)$
depend only on the $p$-bar weight $d$, not on the type $\eps$, of the block.

\section{The affine algebras}\label{s1}

We begin by recalling the construction of the affine Lie algebras 
from \cite[Chapter 8]{Kac}.
Let $X_N^{(r)}$ be an affine Dynkin diagram of ADE type as in the introduction,
and let $X_N$ be the underlying finite Dynkin diagram.
We use the same numbering of Dynkin diagrams as \cite[$\S$4.8]{Kac}
with two exceptions:
in the case $X_N^{(r)} = E_6^{(2)}$ we will number the vertices of the
finite Dynkin diagram $X_N = E_6$ by
$$
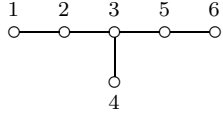
{\begin{picture}(81, 50)%
\put(5, 35){\circle{4}}%
\put(7, 35){\line(1, 0){15}}%
\put(24, 35){\circle{4}}%
\put(26, 35){\line(1, 0){15}}%
\put(43, 35){\circle{4}}%
\put(45, 35){\line(1, 0){15}}%
\put(62, 35){\circle{4}}%
\put(64, 35){\line(1, 0){15}}%
\put(81, 35){\circle{4}}%
\put(43, 33){\line(0, -11){15}}%
\put(43, 16){\circle{4}}%
\put(5, 41){\makebox(0, 0)[b]{$_1$}}%
\put(24, 41){\makebox(0, 0)[b]{$_2$}}%
\put(43, 41){\makebox(0, 0)[b]{$_3$}}%
\put(62, 41){\makebox(0, 0)[b]{$_5$}}%
\put(81, 41){\makebox(0, 0)[b]{$_6$}}%
\put(43, 11){\makebox(0, 0)[t]{$_4$}}%
\end{picture}}
$$
and in the case 
$X_N^{(r)} = A_{2\ell}^{(2)}$ we will number the vertices
of the finite Dynkin diagram $X_N = A_{2\ell}$
by
$$
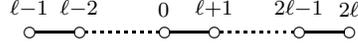
{\begin{picture}(126, 25)%
\put(5, 5){\circle{4}}%
\put(7, 5){\line(1, 0){15}}%
\put(24, 5){\circle{4}}%
\put(26, 5){\line(1, 0){1}}%
\put(29, 5){\line(1, 0){1}}%
\put(32, 5){\line(1, 0){1}}%
\put(35, 5){\line(1, 0){1}}%
\put(38, 5){\line(1, 0){1}}%
\put(41, 5){\line(1, 0){1}}%
\put(44, 5){\line(1, 0){1}}%
\put(47, 5){\line(1, 0){1}}%
\put(50, 5){\line(1, 0){1}}%
\put(53, 5){\line(1, 0){1}}%
\put(56, 5){\circle{4}}%
\put(58, 5){\line(1, 0){15}}%
\put(75, 5){\circle{4}}%
\put(77, 5){\line(1, 0){1}}%
\put(80, 5){\line(1, 0){1}}%
\put(83, 5){\line(1, 0){1}}%
\put(86, 5){\line(1, 0){1}}%
\put(89, 5){\line(1, 0){1}}%
\put(92, 5){\line(1, 0){1}}%
\put(95, 5){\line(1, 0){1}}%
\put(98, 5){\line(1, 0){1}}%
\put(101, 5){\line(1, 0){1}}%
\put(104, 5){\line(1, 0){1}}%
\put(107, 5){\circle{4}}%
\put(109, 5){\line(1, 0){15}}%
\put(126, 5){\circle{4}}%
\put(5, 11){\makebox(0, 0)[b]{$_{\ell-1}$}}%
\put(24, 11){\makebox(0, 0)[b]{$_{\ell - 2}$}}%
\put(56, 11){\makebox(0, 0)[b]{$_{0}$}}%
\put(75, 11){\makebox(0, 0)[b]{$_{\ell+1}$}}%
\put(107, 11){\makebox(0, 0)[b]{$_{2\ell-1}$}}%
\put(126, 11){\makebox(0, 0)[b]{$_{2\ell}$}}%
\end{picture}}
$$

Let $Q'$ denote the root lattice of type $X_N$, with simple roots
$\alpha_i'$ and invariant
bilinear form $(.|.)'$ normalized so that each
$(\alpha_i'| \alpha_i')' = 2$.
Let $\mu:Q' \rightarrow Q'$ be a graph automorphism of order $r$,
as in e.g. \cite[$\S$7.9]{Kac}.
Let 
$$
\eps:Q' \times Q' \rightarrow \{\pm 1\}
$$ 
be an asymmetry function
as in \cite[$\S$7.8]{Kac} chosen so that 
$\eps(\mu(\alpha'), \mu(\beta')) = \eps(\alpha',\beta')$.
In case $X_N^{(r)} = A_{2\ell}^{(2)}$ this is not possible so we instead require here that
$\eps(\mu(\alpha'), \mu(\beta')) = \eps(\beta',\alpha')$.
Let $\mathfrak h' = \C \otimes_\Z Q'$ viewed as an abelian Lie algebra,
and extend $\mu$ and $(.|.)'$ linearly to $\mathfrak h'$.
Then we can construct the finite dimensional simple Lie algebra $\mathfrak g'$
of type $X_N$ as the vector space
$$
\mathfrak g' = \mathfrak h' \oplus \bigoplus_{\alpha' \text{ a root}} \C 
E_{\alpha'}
$$
viewed as a Lie algebra so that
$\mathfrak h'$ is abelian and
\begin{align*}
[\alpha', E_{\beta'}] &= (\alpha'|\beta')' E_{\beta'},
\qquad[E_{\alpha'}, E_{-\alpha'}] = -\alpha',\\
[E_{\alpha'},E_{\beta'}] &= \left\{
\begin{array}{ll}
\eps(\alpha',\beta') E_{\alpha'+\beta'}&\hbox{if $\alpha'+\beta'$ is a root,}\\
0&\hbox{otherwise.}
\end{array}\right.
\end{align*}
The invariant form on $\mathfrak h'$ extends to $\mathfrak g'$
by $(\mathfrak h'|E_{\alpha'})' = 0$ and $(E_{\alpha'}|E_{\beta'})' = 
- \delta_{\alpha',-\beta'}$ for all roots $\alpha',\beta'$.

Let $a_i, a_i^\vee\:(i = 0,\dots,\ell)$ be the numerical labels 
on the
Dynkin diagram $X_N^{(r)}$ and its dual as in \cite[$\S$4.8]{Kac}.
We note especially that $a_0 = 1$ if
$X_N^{(r)} \neq A_{2\ell}^{(2)}$ and
$a_0 = 2$ if $X_N^{(r)} = A_{2\ell}^{(2)}$.
It will also be convenient to define 
\begin{align*}
c_i &= 
\left\{
\begin{array}{ll}
2&\hbox{if $X_N^{(r)} = A_{2\ell}^{(2)}$ and $i = 0$,}\\
1&\hbox{otherwise;}
\end{array}
\right.\\
d_i &= c_i a_i^\vee a_i^{-1} \in \{1,r\}
\end{align*}
for $i = 0,1,\dots,\ell$.
Let $m = a_0 r$ and fix a primitive $m$th root of unity
$\omega \in \C$. 
In all types other than $A_{2\ell}^{(2)}$, 
let $\eta:Q' \rightarrow \C^\times$ denote the constant function with
$\eta(\alpha') = 1$ for all $\alpha' \in Q'$;
in type $A_{2\ell}^{(2)}$, define $\eta$ instead
by the rules
$$
\eta(0) = 1,
\quad
\eta(\alpha'+\beta') = \eta(\alpha') \eta(\beta') (-1)^{(\alpha'|\beta')'},
\qquad
\eta(\alpha_j') = \left\{
\begin{array}{ll}
1&j \neq 0, \ell+1,\\
\omega
&
j = 0, \ell+1.
\end{array}\right.
$$
Now extend $\mu$ from $\mathfrak h'$ to $\mathfrak g'$
by declaring that
$\mu(E_{\alpha'}) = \eta(\alpha') E_{\mu(\alpha')}$ for all roots $\alpha' 
\in Q'$.
The order of the resulting
automorphism $\mu$ of $\mathfrak g'$ is equal to $m$ in all cases.

Decompose 
$$
\mathfrak g' = \bigoplus_{n \in \Z / m} \mathfrak g'_n
\qquad\hbox{where}\qquad
\mathfrak g'_n = \{X \in \mathfrak g'\:|\:\mu(X) = \omega^n X\}.
$$
Also write
$\mathfrak h'_n = \mathfrak h' \cap \mathfrak g'_n$.
Introduce the infinite dimensional Lie algebras
\begin{align*}
{\mathfrak g} &= \bigoplus_{n \in \Z} \mathfrak g'_n \otimes t^n
\oplus \C c \oplus \C d \subseteq \mathfrak g' \otimes \C[t,t^{-1}] \oplus \C c \oplus \C d,\\
{\mathfrak h} &= \mathfrak h'_0 \otimes 1
\oplus \C c \oplus \C d \subset {\mathfrak g},\\
\mathfrak t &= \mathfrak t^+ \oplus \C c \oplus \mathfrak t^- \subset {\mathfrak g}
\quad\hbox{where}\quad
\mathfrak t^{\pm} = \bigoplus_{\pm n > 0} \mathfrak h'_n \otimes t^n.
\end{align*}
Multiplication is defined by the rules
\begin{align*}
& [d, X\otimes t^n ] = n X \otimes t^n,
\qquad
[c, {\mathfrak g}] = 0,\\
& [X\otimes t^n , Y \otimes t^k ] =
[X,Y]\otimes t^{n+k} + \delta_{n,-k} n\frac{(X|Y)'}{m} c.
\end{align*}
Then ${\mathfrak{g}}$ is the affine Lie algebra of type
$X_N^{(r)}$ with canonical central element $c$ and scaling element $d$,
and
${\mathfrak{h}}$ is a Cartan subalgebra.
As a matter of notation, we will write
\begin{equation*}\label{xndef}
X(n) := 
\sum_{j=0}^{m-1} \omega^{-nj} \mu^j(X) \otimes t^n \in \mathfrak g'_n \otimes t^n
\end{equation*}
for $X \in \mathfrak g'$ and $n \in \Z$.
The normalized invariant form on $\mathfrak g$ will be denoted
$(.|.)$, and is defined by
$$
(X \otimes t^n| Y \otimes t^k)
= 
\delta_{n,-k} (X|Y)' / r
$$
for all $X \in \mathfrak g_n, Y \in \mathfrak g_k$.

In order to write down a choice of Chevalley generators for
${\mathfrak{g}}$,
let $\ell$ denote the number of $\mu$-orbits on the simple roots in $Q'$.
Let 
$$
\eps = 
\left\{
\begin{array}{ll}
0 &\hbox{if $X_N^{(r)} \neq A_{2\ell}^{(2)}$},\\
\ell &\hbox{if $X_N^{(r)} = A_{2\ell}^{(2)}$},
\end{array}
\right.
$$ 
and set
$$
I = \{0,1,\dots,\ell\} - \{\eps\}.
$$
Then, the $\alpha_i'$ for
$i \in I$ give a set
of representatives for the $\mu$-orbits on the simple roots.
Define
$$
-\alpha_\eps' = \left\{
\begin{array}{ll}
\hbox{the longest root in $Q'$}&\hbox{if $r = 1$ or $X_N^{(r)} = A_{2\ell}^{(2)}$,}\\
\alpha_1'+\dots+\alpha_{2\ell-2}'&\hbox{if $X_N^{(r)}  = A_{2\ell-1}^{(2)}$,}\\
\alpha_1'+\dots+\alpha_\ell'&\hbox{if $X_N^{(r)}  = D_{\ell+1}^{(2)}$,}\\
\alpha_2'+\alpha_3'+\alpha_4'&\hbox{if $X_N^{(r)} = D_4^{(3)}$,}\\
\alpha_1'+2\alpha_2' + 2\alpha_3' + \alpha_4' + \alpha_5' + \alpha_6'
&\hbox{if $X_N^{(r)} = E_6^{(2)}$}.\\
\end{array}
\right.
$$
For $i = 0,1,\dots,\ell$, write
$$
e_i(n) = \frac{\sqrt{c_i}}{a_0 d_i} E_{\alpha_i'}(n)
\quad\hbox{and}\quad
f_i(n) = -\frac{\sqrt{c_i}}{a_0 d_i} E_{-\alpha_i'}(n).
$$
The Chevalley generators of ${\mathfrak g}$ are
$e_0 = e_0(1), e_i = e_i(0)$ and $f_0 = f_0(-1), f_i = f_i(0)$
for $i = 1,\dots,\ell$, as is proved in \cite[$\S$8.7]{Kac}
(taking $s_0 = 1, s_1=\dots=s_\ell = 0$).
We also define 
$$
h_i = [e_i, f_i] = 
\delta_{i,0} c + \frac{c_i}{a_0 d_i} \alpha_i'(0).
$$

Next let $Q 
\subset {\mathfrak h}^*$ denote the root lattice associated
to ${\mathfrak g}$. So following \cite[$\S$6.2]{Kac},
$$
Q = \bigoplus_{i=0}^\ell \Z {\alpha}_i 
\oplus \Z \Lambda_0
$$ where
$\alpha_0,\dots,\alpha_\ell$ are the simple roots
corresponding to $h_0,\dots,h_\ell$ and $\Lambda_0$ is the zeroth fundamental
dominant weight, i.e.
\begin{align*}
\langle h_i, \alpha_j \rangle 
&= \hbox{the $ij$-entry of the Cartan matrix of type $X_{N}^{(r)}$},\\
\langle h_i, \Lambda_0 \rangle &= 
\langle d, \alpha_i \rangle = \delta_{i,0},\\
\langle d, \Lambda_0 \rangle &= 0,
\end{align*}
for $i,j = 0,\dots,\ell$.
Also as in \cite[$\S$6.2]{Kac}, we 
have the normalized invariant form $(.|.)$ on $\mathfrak h^*$
and the element $\delta = \sum_{i=0}^\ell a_i \alpha_i \in Q$.

To conclude, we explain the relationship between the
form $(.|.)'$ on $Q'$ and the form $(.|.)$ on $Q$.
Introduce the new symmetric bilinear form $(.|.)_\mu$ on $Q'$ 
defined by
\begin{equation*}\label{newform}
(\alpha'|\beta')_\mu = (\alpha'|\sum_{j=0}^{r-1} \mu^j(\beta'))'
\end{equation*}
for all $\alpha',\beta' \in Q'$. 
There is an orthogonal decomposition
$$
\mathfrak h^* = \ocirc{\mathfrak h}{\!\!^*} \oplus
(\C \delta + \C \Lambda_0)
$$
where $\ocirc{\mathfrak h}{\!\!^*} = \bigoplus_{i=1}^\ell \C \alpha_i$,
see \cite[$\S$6.2]{Kac}.
As in {\em loc. cit.} we write $-:\mathfrak h^* \rightarrow
\ocirc{\mathfrak h}{\!\!^*}$ for the orthogonal projection, in particular
$\overline{Q}$ denotes the orthogonal projection of $Q$ onto
$\ocirc{\mathfrak h}{\!\!^*}$.
Define a $\Z$-linear map 
\begin{equation}\label{atend}
\iota:Q' \rightarrow \overline{Q}
\end{equation}
by 
$\iota(\mu^j(\alpha_i')) = 
\overline{\alpha}_i$
for each $i \in I$ and $j \geq 0$.
The kernel of $\iota$ is the space
\begin{equation}
\label{EM}
M' = \{\alpha' - \mu(\alpha')\:|\:\alpha' \in Q'\}
\end{equation}
which is precisely the radical of the bilinear form $(.|.)_\mu$.
Moreover, $\iota$ induces an isometry between $Q' / M'$
and $\overline{Q}$
with respect to the forms induced by $(.|.)_\mu$ and $(.|.)$ 
respectively.

\section{The basic representation}

Next we recall the construction of the basic representation $V = V(\Lambda_0)$
of $\mathfrak g$, following Lepowsky \cite{Lep}.
Let $Z = \langle -1, \omega \rangle \subset \C^\times$ be the multiplicative group generated by $-1$ and $\omega$.
Form the central extension
$$
1 \longrightarrow Z \longrightarrow \widehat Q \stackrel{\pi}{\longrightarrow} 
Q' \longrightarrow 1,
$$
namely,
$\widehat Q = \{e^{\alpha'}_x\:|\:\alpha' \in Q', x \in Z\}$
with multiplication
$$
e^{\alpha'}_x e^{\beta'}_y = 
\left\{
\begin{array}{ll}
e^{\alpha'+\beta'}_{xy \eps(\alpha',\beta')}&\hbox{if $X_N^{(r)} \neq A_{2\ell}^{(2)}$, $D_4^{(3)}$,}\\
e^{\alpha'+\beta'}_{xy \eps(\alpha',\beta')(-\omega)^{-(\alpha'|\mu(\beta'))'}}
&\hbox{if $X_N^{(r)} = A_{2\ell}^{(2)}$ or $D_4^{(3)}$,}
\end{array}\right.
$$
for $\alpha',\beta' \in Q',\ x,y \in Z$.
The map $\pi:\widehat Q \rightarrow Q'$ here 
is defined by $\pi(e^{\alpha'}_x) = \alpha'$.
Let $\widehat M = \pi^{-1}(M')$, where $M'$ is as in (\ref{EM}).
There is a well-defined 
multiplicative character $\tau:\widehat M \rightarrow \C^\times$
defined in \cite[Proposition~6.1]{Lep} by
$$
\tau(e^{\alpha' - \mu(\alpha')}_x)
= (-1)^{(\alpha'|\alpha')'/2} x \eta(\alpha') \eps(\alpha',\mu(\alpha'))
\omega^{-a_0^2(\alpha'|\alpha')_\mu/2}.
$$
So we can form the induced $\widehat Q$-module
$\C [\widehat Q] \otimes_{\C [\widehat M]} \tau.$
We note the useful formula 
\begin{equation*}
\label{EUseful}
e^{\alpha'}_1 \otimes \tau = 
\eta(\alpha') \omega^{a_0(\alpha'|\alpha')_\mu/2}
e^{\mu(\alpha')}_1 \otimes \tau \qquad (\alpha' \in Q').
\end{equation*}

View the symmetric algebra $S(\mathfrak t^-)$ as a $\mathfrak t$-module
in the unique way so that $c$ acts as $1$, elements of $\mathfrak t^-$ act 
by multiplication, and elements of $\mathfrak t^+$ annihilate $1$.
It is $\Z$-graded by declaring that $$
\deg(h \otimes t^{-n}) = \frac{n}{a_0}
$$
for each $h \in \mathfrak h_{-n}', n \geq 1$.
Let
$$
V = S(\mathfrak t^-) \otimes \C [\widehat Q] \otimes_{\C [\widehat M]} \tau.
$$
Let $\mathfrak t$ 
act on $S(\mathfrak t^-)$ as given 
and trivially on $\C[\widehat Q] \otimes_{\C[\widehat M]} \tau$, 
let $h \otimes t^0$ for $h \in \mathfrak h_0'$ act by
$$
(h \otimes t^0) (f \otimes e_x^{\alpha'} \otimes \tau) = (h|\alpha')' 
f \otimes e_x^{\alpha'}\otimes \tau,
$$
and let $d$ act by
$$
d(f \otimes e_x^{\alpha'} \otimes \tau) = -a_0\left(\deg(f) 
+ (\alpha'|\alpha')_\mu/2\right)f \otimes e_x^{\alpha'} \otimes \tau.
$$
We have now defined the action of ${\mathfrak h} + \mathfrak t$ on $V$.
To extend the action to all of ${\mathfrak g}$, 
let $\alpha' \in Q'$ be a root.
As in \cite[(4.8)]{Lep}, let
$$
\sigma(\alpha')
=
\left\{
\begin{array}{ll}
1&\hbox{$r=1$,}\\
\sqrt{2}^{(\alpha'|\mu(\alpha'))'}&\hbox{if $X_N^{(r)} = A_{2\ell-1}^{(2)}, 
D_{\ell+1}^{(2)}$ or $E_6^{(2)}$,}\\
(1-\omega^{-1})^{(\alpha'|\mu(\alpha'))'}&\hbox{if $X_N^{(r)} = D_4^{(3)}$,}\\
2(1+\omega)^{(\alpha'|\mu(\alpha'))'}&\hbox{if $X_N^{(r)} = A_{2\ell}^{(2)}$}.
\end{array}
\right.
$$
Also define
$$
P_{\alpha'}(z) = 
\exp\left(\sum_{n \geq 1} \frac{\alpha'(- n) z^n}{n}\right),
\qquad
Q_{\alpha'}(z) = 
\exp\left(-\sum_{n \geq 1} \frac{\alpha'(n) z^n}{n}\right),
$$
viewed as elements of $\End(V)[[z^{\pm 1}]]$.
Let
$$
E_{\alpha'}(z) = \sigma(\alpha')P_{\alpha'}(z) Q_{\alpha'}(z^{-1})
e_1^{\alpha'}
z^{a_0\alpha'} z^{a_0(\alpha'|\alpha')_\mu/2-1}.
$$
Here, $z^{a_0\alpha'}$ denotes the operator with
$$
z^{a_0\alpha'} (f \otimes e_x^{\beta'} \otimes \tau) 
= z^{(a_0\alpha'|\beta')_\mu} f \otimes e_x^{\beta'} \otimes \tau
$$ 
for each $f \in S(\mathfrak t^-)$ and $\beta \in Q'$,  and
$$
e_1^{\alpha'}(f \otimes e_x^{\beta'}\otimes\tau) =
f\otimes (e_1^{\alpha'} e_x^{\beta'}) \otimes \tau.
$$
Expanding $E_{\alpha'}(z)$ in powers of $z$ we get the required action of
$E_{\alpha'}(n) 
\in {\mathfrak g}$  on $V$ for each root $\alpha' \in Q'$ 
and each $n \in \Z$:
$$
E_{\alpha'}(z) = \sum_{n \in \Z} E_{\alpha'}(n) z^{-n-1}.
$$
For a proof that this is a well-defined irreducible representation of
${\mathfrak g}$ in case $r = 1$ see \cite[$\S$14.8]{Kac};
the general case is due to Lepowsky \cite{Lep}.

Let $\C[\overline Q]$ denote the group algebra of $\overline Q$, 
with natural basis
$e^{\alpha}$ for $\alpha \in \overline{Q}$
and multiplication
$e^{\alpha} e^{\beta} = e^{\alpha+\beta}$.
Note $\C[\widehat Q] \otimes_{\C[\widehat M]} \tau$ has
a basis given by the elements $e^{\alpha'}_1 \otimes \tau$ for
all $\alpha' \in \sum_{i \in I} \Z \alpha_i'$.
For such an $\alpha'$, let
$$
\iota(e_1^{\alpha'} \otimes \tau) = 
\left\{
\begin{array}{ll}
e^{\iota(\alpha')}&\hbox{if $X_N^{(r)} \neq A_{2\ell}^{(2)}$, $D_4^{(3)}$,}\\
(-\omega)^{(\alpha'|\mu(\alpha'))'/2}e^{\iota(\alpha')}
&\hbox{if $X_N^{(r)} = D_4^{(3)}$,}\\
\left(\frac{1-\omega}{\sqrt 2} \right)^{(\alpha'|\mu(\alpha'))'}
e^{\iota(\alpha')}
&\hbox{if $X_N^{(r)} = A_{2\ell}^{(2)}$,}
\end{array}\right.
$$
recalling the map $\iota:Q' \rightarrow \overline{Q}$ defined in (\ref{atend}).
Extending linearly, we obtain a vector space isomorphism
$\iota:\C [\widehat Q] \otimes_{\C [\widehat M]} \tau
 \rightarrow \C[\overline{Q}]$.
For $i = 0,1,\dots,\ell$, we 
define functions $\sigma_i^{\pm}:\overline{Q} \rightarrow \C^\times$ 
by the equation
\begin{equation*}\label{sdef}
\sigma^{\pm}_i(\alpha) e^{\alpha \pm \overline{\alpha}_i}
=
\pm \frac{\sqrt{c_i}}{a_0 d_i}\sigma(\alpha_i') 
\iota(e_1^{\pm \alpha_i'} \iota^{-1}(e^{\alpha}))
\end{equation*}
for all $\alpha \in \overline{Q}$.
The choice of the renormalization map $\iota$ above ensures:

\vspace{1mm}
\begin{Lemma}\label{sigmal}
For all $i = 0,1,\dots,\ell$ and $\alpha \in \overline{Q}$, 
$\sigma_i^{\pm}(\alpha) \in \{\pm 1\}$.
Moreover, for $i \in I$, we have that $\sigma_i^- = - \sigma_i^+$,
and $\sigma_i^+:\overline{Q} \rightarrow \{\pm 1\}$
is a group homomorphism such that
$\sigma_i^+(\overline{\alpha}_j) = \eps(\alpha_i', \alpha_j')$
for each $j \in I$.
\end{Lemma}
\vspace{1mm}

Now we can rewrite the construction of the basic representation $V$
in terms of the Chevalley generators.
We will identify 
$$
V = S(\mathfrak t^-) \otimes \C[\widehat{Q}]\otimes_{\C[\widehat{M}]} \tau 
= S(\mathfrak t^-) \otimes \C[\overline{Q}]
$$
via the map $\id\otimes\iota$.
Then,
the actions of $h_i$ for $i = 0,\dots,\ell$ and of $d$ are as
\begin{align*}
h_i (f \otimes e^{\alpha}) &= (\delta_{i,0} + \langle h_i, \alpha \rangle) f \otimes e^{\alpha},\\
d(f \otimes e^{\alpha})&= - a_0\left(\deg(f) +  
(\alpha|\alpha)/2\right) f \otimes e^{\alpha}
\end{align*}
for all $\alpha \in \overline{Q}$.
In particular, we note from this that
\begin{equation}
\label{EWeight}
\wt(f \otimes e^{\alpha}) = 
\Lambda_0 + \alpha  - \left(\deg(f) +
(\alpha|\alpha)/2
\right)\delta
\end{equation}
for each homogeneous $f \in S(\mathfrak t^-)$ and $\alpha \in \overline{Q}$. 
This shows that $1 \otimes e^{0}$ is a highest weight vector
in $V$ of highest weight $\Lambda_0$ (cf. \cite[Lemma 12.6]{Kac}), 
identifying $V$ with the irreducible
highest weight module $V(\Lambda_0)$.
Finally, for $i = 0,\dots,\ell$,
\begin{align}
e_i(z) &= \sum_{n \in \Z} e_i(n) \otimes z^{-n-1} =
P_{\alpha_i'}(z) Q_{\alpha_i'}(z^{-1}) e^{\overline{\alpha}_i} z^{a_0 {\alpha}_i} 
z^{{a_0}({\alpha}_i|{\alpha}_i)/2 - 1} s^+_i,\label{v1}\\\label{v2}
f_i(z) &= \sum_{n \in \Z} f_i(n) \otimes z^{-n-1} =
P_{-\alpha_i'}(z) Q_{-\alpha_i'}(z^{-1}) e^{-\overline{\alpha}_i} z^{-a_0{\alpha}_i} 
z^{{a_0}({\alpha}_i|{\alpha}_i)/2 - 1} s^-_i,
\end{align}
where
\begin{align*}
z^{\pm a_0{\alpha}_i} (f \otimes e^{\beta})
&= z^{(\pm a_0{\alpha}_i|\beta)} f \otimes e^{\beta},\\
s^{\pm}_i(f \otimes e^{\beta}) &= \sigma^{\pm}_i(\beta) f \otimes e^{\beta},\\
e^{\pm\overline{\alpha}_i}(f \otimes e^{\beta}) &= 
f \otimes e^{\beta \pm \overline{\alpha}_i}.
\end{align*}
The following lemma will be needed later on:

\vspace{1mm}
\begin{Lemma}\label{bigprod}
For $i_1,\dots,i_s \in \{0,\dots,\ell\}$,
roots $\beta_1',\dots,\beta_t' \in Q'$
and 
$\gamma\in \overline{Q}$, we have that
\begin{eqnarray*}
& & 
\hspace {-1.3 cm}e_{i_1}(z_1)e_{i_2}(z_2)
\dots e_{i_s}(z_s) P_{\beta_1'}(w_1) \dots P_{\beta_{t}'}(w_t)
\otimes e^\gamma = \\
& & 
\pm 
\prod_{1\leq u \leq s} 
z_u^{\frac{a_0}{2}(\alpha_{i_u}|\alpha_{i_u}) - 1 +
a_0 (\alpha_{i_u}|\gamma)}
\\
& &
\times 
\prod_{1 \leq u  < v \leq s}
z_u^{a_0(\alpha_{i_u}|\alpha_{i_v})}
\prod_{k\in \Z/m}\left(1 - \omega^{-k} 
\frac{z_v}{z_u}\right)^{(\mu^k(\alpha_{i_u}')|\alpha_{i_v}')'}
\\
& &
\times 
\prod_{1\leq u \leq s}\prod_{1 \leq v \leq t}
\prod_{k \in \Z / m}
\left(1 - \omega^{-k} \frac{w_v}{z_u}\right)^{(\mu^k(\alpha_{i_u}')|\beta_v')'}
\\
& &
\times P_{\alpha_{i_1}'}(z_1) \dots P_{\alpha_{i_s}'}(z_s)
P_{\beta_1'}(w_1) \dots P_{\beta_t'}(w_t) \otimes 
e^{\gamma+\overline{\alpha}_{i_1} + \dots +\overline{\alpha}_{i_s}}.\quad
\end{eqnarray*}
A similar formula holds for $f_{i_1}(z_1) \dots f_{i_s}(z_s)
P_{\beta_1'}(w_1) \dots P_{\beta_{s}'}(w_s)
\otimes e^\gamma$, replacing 
$\alpha_{i_u}$ by $- \alpha_{i_u}$,
$\alpha_{i_u}'$ by $- \alpha_{i_u}'$, and $\bar\alpha_{i_u}$ by $- \bar\alpha_{i_u}$
everywhere.
\end{Lemma}

\begin{proof}
This follows from the following commutation relation
obtained in \cite[3.4]{Lep}: for $\alpha',\beta' \in Q'$,
\begin{equation*}\label{useful}
Q_{\alpha'}(z^{-1}) P_{\beta'}(w)
= 
P_{\beta'}(w)Q_{\alpha'}(z^{-1}) 
\prod_{k\in \Z / m} 
(1 - \omega^{-k} \frac{w}{z})^{(\mu^k(\alpha')|\beta')'},
\end{equation*}
which is a consequence of the Campbell-Hausdorff formula,
cf. \cite[(14.8.12)]{Kac}.
\end{proof}

\section{The Integral form}

As in the introduction, 
let $U_\Z$ denote the $\Z$-subalgebra of the universal enveloping
algebra of ${\mathfrak g}$ generated by the elements
$e_i^r / r!, f_i^r / r!$ for $i = 0,\dots,\ell$ and $r \geq 0$, and
let 
$$
V_\Z := U_\Z (1 \otimes e^0) \subset V.
$$ 
In this section, we will give an explicit description
of $V_\Z$.

To start with,
let $\tau:\mathfrak g \rightarrow \mathfrak g$ be the antilinear Chevalley 
antiautomorphism, so 
$$
\tau(d) = d,\quad 
\tau(e_i(n)) = f_i(-n),\quad \tau(f_i(n)) = e_i(-n)
$$
for each $i = 0,\dots,\ell$ and $n \in \Z$, cf. \cite[$\S\S$7.6, 8.3]{Kac}. 
The {\em Shapovalov form}\, $(.|.)_S$ on $V$ is the unique Hermitian form
such that
$$
(1 \otimes e^{0}, 1 \otimes e^{0})_S = 1\quad \text{and} \quad
(x v, w)_S = (v, \tau(x) w)_S
$$ 
for all $v,w \in V$, $x \in \mathfrak g$.
The restriction of $\tau$ to $\mathfrak t$ gives the antilinear Chevalley
antiautomorphism of $\mathfrak t$, and we can also 
consider the Shapovalov form on $S(\mathfrak t^-)$,
satisfying
$(1,1)_S = 1$ and $(x f, g)_S = (f, \tau(x) g)_S$
for all $f,g \in S(\mathfrak t^-)$, $x \in \mathfrak t$.

\vspace{1mm}
\begin{Lemma}\label{one} For all $f,g \in S(\mathfrak t^-)$ and
$\alpha,\beta \in \overline{Q}$,
$(f \otimes e^{\alpha}, g \otimes e^{\beta})_S
= (f,g)_S.$
\end{Lemma}

\begin{proof}
Since different weight spaces are orthogonal and in view of (\ref{EWeight}), 
this reduces to checking that
$(1 \otimes e^{\alpha}, 1 \otimes
e^{\alpha})_S = 1$ for all $\alpha \in \overline{Q}$.
Proceeding by induction, we may assume that
there is some $\beta \in \overline{Q}$ and $i \in I$ such that
$(1 \otimes e^{\beta},1 \otimes e^{\beta})_S = 1$
and 
either $\alpha = \beta + \overline{\alpha}_i$
or $\alpha = \beta - \overline{\alpha}_i$.

Suppose that $\alpha = \beta + \overline{\alpha}_i$.
Letting $n = - a_0(\alpha_i|\beta) - a_0(\alpha_i|\alpha_i)/2$,
one checks easily using (\ref{v1}), (\ref{v2}) that
\begin{equation}\label{inpf}
e_i(n) (1 \otimes e^{\beta}) = \sigma_i^+(\alpha)( 1 \otimes e^{\alpha}),\qquad
f_i(-n) (1 \otimes e^{\alpha}) = \sigma_i^-(\beta) (1 \otimes e^{\beta}).
\end{equation}
Hence, 
\begin{align*}
(1 \otimes e^{\alpha}, 1 \otimes e^{\alpha})_S
&=
\sigma^+_i(\beta)(e_i(n)(1 \otimes e^{\beta}),1 \otimes e^{\alpha})_S\\
&=
\sigma^+_i(\beta)(1 \otimes e^{\beta}, f_i(-n)(1 \otimes e^{\alpha}))_S\\
&=
\sigma_i^-(\alpha) \sigma_i^+(\beta) (1 \otimes e^{\beta}, 1 \otimes 
e^{\beta})_S = 1,
\end{align*}
since $\sigma_i^-(\alpha) = \sigma_i^+(\beta)$ by Lemma~\ref{sigmal}.

A similar argument in the case that $\alpha = \beta - \overline{\alpha}_i$
completes the proof.
\end{proof}

\begin{Lemma}\label{two}
For all $i =0,1,\dots,\ell$ and $n \in \Z$, the elements
$e_i(n)$ and $f_i(n)$ belong to $U_\Z$.
\end{Lemma}

\begin{proof}
Suppose that $e_i(n) \neq 0$.
Then, $\wt(e_i(n)) = \overline{\alpha}_i + \frac{n}{a_0} \delta$ is a real
root, hence is conjugate under the Weyl group $W$ associated to
$\mathfrak g$ to some simple root $\alpha_j$.
So we can find simple reflections $s_{i_1},\dots,s_{i_t} \in W$
such that $\overline{\alpha}_i + \frac{n}{a_0} \delta
= s_{i_1} \dots s_{i_t} \alpha_j$.
Let $r^\ad_i$ be the automorphism of $\mathfrak g$
defined by
$r_i^\ad = \exp (\ad f_i) \exp (-\ad e_i) \exp (\ad f_i)$, for $i = 0,1,\dots,\ell$.
Since real root spaces of $\mathfrak g$ are one dimensional,
$$
r_{i_1}^\ad \dots r_{i_t}^{\ad} e_j = c e_i(n)
$$
for some non-zero scalar $c$. 
Now, $\tau(\exp(\ad y) (x))=\exp(-\ad\tau(y))(\tau(x))$, whence by an $SL_2$-calculation we have $r_i^\ad (x) = \tau(r_i^\ad(\tau (x)))$ for all $x \in 
\mathfrak g$, we also get that
$$
r_{i_1}^\ad \dots r_{i_t}^{\ad} f_j = c f_i(-n).
$$
But the $r_i^\ad$ preserve the normalized invariant form on $\mathfrak g$,
so
$$
a_j (a_j^{\vee})^{-1} = (e_j| f_j) = (c e_i(n)|c f_i(-n)) = c^2
a_i (a_i^{\vee})^{-1}.
$$
Clearly, $\alpha_i$ and $\alpha_j$ are roots of the same length, i.e.
$a_j (a_j^{\vee})^{-1} = 
a_i (a_i^{\vee})^{-1}$, so this gives that $c = \pm 1$.
Finally, the action of $r^\ad_i$ on $\mathfrak g$ 
leaves $U_\Z \cap \mathfrak g$ invariant, so
$$
e_i(n) = \pm r^\ad_{i_1} \dots r^\ad_{i_s} e_j \in U_\Z,
$$
and similarly $f_i(n) \in U_\Z$ too.
\end{proof}

For $n \geq 1$ and $i = 0,1,\dots,\ell$, define
$$
y_{nd_i}^{(i)} = \frac{\alpha_i'(-a_0n d_i)}{a_0 nd_i},
\qquad
x_{nd_i}^{(i)} = \sum_{k_1+2k_2+\dots = n}
\frac{{y_{d_i}^{(i)}}^{k_1}}{k_1!}
\frac{{y_{2d_i}^{(i)}}^{k_2}}{k_2!}
\frac{{y_{3d_i}^{(i)}}^{k_3}}{k_3!}\cdots.
$$
Observe that
\begin{equation}
\label{EYX}
P_{\alpha_i'}(z) = \exp\left(\sum_{n \geq 1} y_{nd_i}^{(i)} z^{a_0 n d_i}
\right)
= 1 + \sum_{n \geq 1} x_{nd_i}^{(i)} z^{a_0 n d_i}.
\end{equation}
The $y_{nd_i}^{(i)}$ for $n \geq 1$ and $i \in I$
give a basis
for $\mathfrak t^-$. So
$S(\mathfrak t^-)$ is equal to the
free polynomial algebra 
$$
B := \C[y_{nd_i}^{(i)}\:|\:n \geq 1, i \in I].
$$
Since the $x_{nd_i}^{(i)}$ are related to the $y_{nd_i}^{(i)}$ 
in a unitriangular way,
we obtain a $\Z$-form
$$
B_\Z := \Z[x_{nd_i}^{(i)}\:|\:n \geq 1, i \in I]
\subset B
$$
for $B$.
As $\al_\eps'$ is an integral linear combination of the $\al_i'$ with $i\in I$,
it follows from (\ref{EYX}) that the elements $x_{nd_\eps}^{(\eps)}$
also belong to the lattice $B_\Z$.
The $\Z$-grading on $B_\Z$ induced by the grading on $S(\mathfrak t^-)$
is determined by $\deg(y_{n}^{(i)}) = \deg(x_{n}^{(i)}) = n$.

The following theorem (or rather its $q$-analogue) for the non-twisted 
case has been proved in \cite{CJ}. Our argument for the general case is 
similar.

\vspace{1mm}
\begin{Theorem}\label{zform}
$V_\Z = B_\Z \otimes_{\Z} \Z[\overline{Q}]$.
\end{Theorem}

\begin{proof}
Let us first show that
$B_\Z \otimes_{\Z} \Z[\overline{Q}] \subseteq V_\Z$.
Fix $i_1,\dots,i_s \in I$, and let
$$
M(i_1,\dots,i_s) = \{(n_1,\dots,n_s)\:|\:n_1 \geq \dots \geq n_s \geq 0
\hbox{ and }
d_{i_j} | n_j\hbox{ for all $j=1,\dots,s$}\}.
$$
Denote by $>$ the dominance ordering on partitions belonging
to $M(i_1,\dots,i_s)$.
We will show that
$x_{n_1}^{(i_1)} \dots x_{n_s}^{(i_s)} \otimes e^\beta \in V_\Z$
for all $(n_1,\dots,n_s) \in M(i_1,\dots,i_s)$ and each
$\beta \in \overline{Q}$.
Clearly every monomial in $B_\Z$ is of the form
$x_{n_1}^{(i_1)} \dots x_{n_s}^{(i_s)}$ for some choice
of $i_1,\dots,i_s$ and $(n_1,\dots,n_s) \in M(i_1,\dots,i_s)$,
so this is good enough.

To start with, each $e_i(n), f_i(-n) \in U_\Z$ by Lemma~\ref{two}.
So an obvious inductive argument using (\ref{inpf})
gives that $1 \otimes e^\gamma \in V_\Z$
for each $\gamma \in \overline{Q}$.
Hence, letting $\gamma = \beta - \bar\alpha_{i_1} - \dots - \bar\alpha_{i_s}$,
Lemma~\ref{two} implies that 
all coefficients of $e_{i_1}(z_1) \dots e_{i_s}(z_s) \otimes e^\gamma$
belong to $V_\Z$.
Applying Lemma~\ref{bigprod}, we deduce that all the coefficients of
$$
X := \left(\prod_{1 \leq u  < v \leq s}
\prod_{k\in \Z/m}
\left(1 - \omega^{-k} \frac{z_v}{z_u}\right)^{(\mu^k(\alpha_{i_u}')|\alpha_{i_v}')'}
\right)
P_{\alpha_{i_1}'}(z_1) \dots P_{\alpha_{i_s}'}(z_s) \otimes 
e^{\beta}
$$
belong to $V_\Z$.
One checks that in all cases,
$$
\prod_{k\in \Z/m}
\left(1 - \omega^{-k} \frac{z_v}{z_u}\right)^{(\mu^k(\alpha_{i_u}')|\alpha_{i_v}')'}
= 1 + (*),
$$
where ($*$) is a $\Z$-linear combination of $\left(\frac{z_v}{z_u} \right)^p$
for $p \geq 1$.
It follows that the $z_1^{a_0n_1} \dots z_s^{a_0 n_s}$-coefficient
of $X$ equals
$$
x_{n_1}^{(i_1)} \dots x_{n_s}^{(i_s)}\otimes e^\beta + (**),
$$
where ($**$) is a $\Z$-linear combination
of $x_{n_1'}^{(i_1)} \dots x_{n_s'}^{(i_s)}\otimes e^\beta$ for
$(n_1',\dots,n_s') > (n_1,\dots,n_s)$.
Using downward induction on this ordering,
we deduce that $x_{n_1}^{(i_1)} \dots x_{n_s}^{(i_s)} \in V_\Z$.

Finally, we prove that $B_\Z \otimes_{\Z} \Z[\overline{Q}] \supseteq V_\Z$.
As the high weight vector $1\otimes e^0$ belongs to 
$B_\Z \otimes_{\Z} \Z[\overline{Q}] \supseteq V_\Z$, it 
suffices to show that $B_\Z \otimes_{\Z} \Z[\overline{Q}]$ is invariant 
under each of the
operators $f_i(n)^s / s!$ for $n \in \Z, s \geq 1$
and $i = 0,1,\dots,\ell$. Fix $i\in\{0,1,\dots,\ell\}$
and consider
$$
Y := f_i(z_1)\dots f_i(z_s) P_{\alpha_{i_1}'}(w_1) \dots P_{\alpha_{i_t}'}(w_t)
\otimes e^{\gamma+s \overline{\alpha}_i}
$$
for $i_1,\dots,i_t \in I$ and $\gamma \in \overline{Q}$.
Applying Lemma~\ref{bigprod} and simplifying, 
\begin{eqnarray*}
Y & = & 
\pm (z_1 \dots z_s)^{\frac{a_0}{2}(\alpha_{i}|\alpha_{i}) - 1 -
a_0 (\alpha_{i}|\gamma + s \overline{\alpha}_i)} \\
& & 
\times 
\prod_{1 \leq u  < v \leq s}
\prod_{k\in \Z/m}
\left(z_u - \omega^{-k} z_v\right)^{(\mu^k(\alpha_{i}')|\alpha_{i}')'}
\\
& & \times 
\prod_{\stackrel{\hbox{$_{1 \leq u \leq s}$}}{1 \leq v \leq t}}
\prod_{k \in \Z / m}
\left(1 - \omega^{-k} \frac{w_v}{z_u}\right)^
{-(\mu^k(\alpha_{i}')|\alpha_{i_v}')'}
\\
& &
\times P_{-\alpha_{i}'}(z_1) \dots P_{-\alpha_{i}'}(z_s)
P_{\alpha_{i_1}'}(w_1) \dots P_{\alpha_{i_t}'}(w_t) \otimes 
e^{\gamma}.
\end{eqnarray*}
Certainly, each coefficient of 
$P_{-\alpha_{i}'}(z_1) \dots P_{-\alpha_{i}'}(z_s)P_{\alpha_{i_1}'}(w_1) \dots P_{\alpha_{i_t}'}(w_t) \otimes 
e^{\gamma}$ belongs to $V_\Z$.
One checks that
$$
\prod_{k\in \Z/m}
\left(z_u - \omega^{-k} z_v\right)^{(\mu^k(\alpha_{i}')|\alpha_{i}')'}
= \left\{
\begin{array}{ll}
(z_u - z_v)^2 \frac{(z_u+z_v)^2}{z_u^2+z_v^2}&\hbox{if $X_N^{(r)} = A_{2\ell}^{(2)}$ and $i = 0$},\\
(z_u^{a_0 d_i} - z_v^{a_0 d_i})^2&\hbox{otherwise.}
\end{array}
\right.
$$
Hence in all cases, $Y$ 
looks like $\prod_{1 \leq u < v \leq s} (z_u - z_v)^2$ times an expression that is symmetric in $z_1,\dots,z_s$.
Hence, by \cite[Lemma 2.5(ii)]{CJ}, the coefficient of $(z_1\dots z_s)^{-n-1}$
in $Y$ is divisible by $s!$. Hence,
all coefficients of 
$(f_i(n)^s / s!) P_{\alpha_{i_1}'}(w_1) \dots P_{\alpha_{i_s}'}(w_s)
\otimes e^{\gamma - s \overline{\alpha}_i}$
belong to $V_\Z$, which completes the proof.
\end{proof}

\section{The determinant}

Fix now some $d \geq 0$.
Lemma~\ref{one} and Theorem~\ref{zform} reduce the problem of computing the
determinant of the Shapovalov form on the $(w \Lambda_0 - d \delta)$ weight
space of $V_\Z$ for any $w \in W$ 
to the problem of computing the determinant of
the Shapovalov form on the degree $d$ component of $B_\Z$.
To tackle the latter question, observe
$$
B = \bigotimes_{i \in I}
\C[y_{d_i}^{(i)},y_{2d_i}^{(i)}, \dots],
\qquad
B_\Z = \bigotimes_{i \in I}
\Z[x_{d_i}^{(i)},x_{2d_i}^{(i)}, \dots].
$$
So according to \cite[Corollary 2.1]{CKK}, there is a well-defined Hermitian 
form on $B$ determined by the rules
$$
(1,1)_K = 1,
\qquad
(n y_{nd_i}^{(i)} f, g)_K = (f, \frac{\partial}{\partial y_{nd_i}^{(i)}} g)_K
$$
for all $i \in I, n \geq 1$ and $f,g \in B$.
Moreover, there is a homogeneous 
basis for the lattice $B_\Z$ (given by Schur polynomials)
that is orthonormal with respect to the form $(.\,,.)_K$.
In particular, the determinant of the form $(.\,,.)_K$ on the degree $d$
 component of $B_\Z$ is equal to $1$.
Our strategy will therefore be to relate the Shapovalov form $(.\,,.)_S$ to
the form $(.\,,.)_K$.

Define $I(n) = \{i \in I\:|\:d_i | n\}$.
Introduce the matrices $A^{(n)} = (a_{i,j}^{(n)})_{i,j \in I(n)}$ 
with
$$
a_{i,j}^{(n)} = \frac{1}{d_i} 
(\alpha_i'|\sum_{k=0}^{r-1}\omega^{a_0nk} \mu^k(\alpha_j'))'.
$$
Recall $\alpha,\beta$ from the introduction. One verifies:

\vspace{1mm}
\begin{Lemma}\label{f0}
For any $n \geq 0$, we have 
$
\det A^{(n)} = \left\{
\begin{array}{ll}
\alpha&\hbox{if $r \mid n$,}\\
\beta&\hbox{if $r \nmid n$,}
\end{array}
\right.
$
\end{Lemma}
\vspace{1mm}

The significance of the matrices $A^{(n)}$ is that for $i \in I(n)$, 
the element
$\tau(n y_{n}^{(i)} / d_i) = \alpha_i'(a_0 n) / a_0 d_i
$
acts on $B$ as the operator
$\sum_{j \in I{(n)}} a^{(n)}_{i,j} 
\frac{\partial}{\partial y_{n}^{(j)}}.$
We set
$$
z_n^{(i)} = \sum_{j \in I(n)} a_{i,j}^{(n)} y_n^{(j)}.
$$
For a partition $\lambda = (\lambda_1 \geq \dots \geq \lambda_h > 0)$
we let $I(\lambda) = I(\lambda_1) \times \dots \times I(\lambda_h)$.
Given $\bi = (i_1,\dots,i_h) \in I(\lambda)$, let
$$
x_\lambda^{(\bi)} = x_{\lambda_1}^{(i_1)}
\dots x_{\lambda_h}^{(i_h)},\qquad
y_\lambda^{(\bi)} = y_{\lambda_1}^{(i_1)}
\dots x_{\lambda_h}^{(i_h)},\qquad
z_\lambda^{(\bi)} = z_{\lambda_1}^{(i_1)}
\dots x_{\lambda_h}^{(i_h)},
$$
all elements of $B$ of degree $|\lambda|$.

\vspace{1mm}
\begin{Lemma}\label{f2}
For $\bi \in I(\lambda)$ and any $f \in B$, we have 
$(y_\lambda^{(\bi)}, f)_S
= (z_\lambda^{(\bi)}, f)_K$.
\end{Lemma}

\begin{proof}
Proceed by induction on the number of non-zero parts of $\lambda$, starting 
induction from the obvious fact that $(1, f)_S=(1, f)_K$.
For the induction step, note that for $i \in I(n)$,
\begin{align*}
(y_{n}^{(i)} y_\lambda^{(\bi)}, f)_S
&=
(\frac{d_i}{n} y_\lambda^{(\bi)},
\sum_{j\in I(n)} {a_{i,j}^{(n)}}
\frac{\partial}{\partial y_{n}^{(j)} }f)_S\\
&=
(\frac{d_i}{n} z_\lambda^{(\bi)},
\sum_{j\in I(n)} {a_{i,j}^{(n)}}
\frac{\partial}{\partial y_{n}^{(j)} }f)_K
=
(z_{n}^{(i)} z_\lambda^{(\bi)}, f)_K.
\end{align*}
\end{proof}

Let
$$
\Omega(\lambda) = \{\bi \in I(\lambda)\:|\:
i_j \leq i_{j+1}\text{ whenever }\lambda_j = \lambda_{j+1}\}.
$$
Then
$\{x_{\lambda}^{(\bi)}\:|\:\lambda \in \mathscr P(d), \bi \in \Omega(\lambda)\}$, 
$\{y_{\lambda}^{(\bi)}\:|\:\lambda \in \mathscr P(d), \bi \in 
\Omega(\lambda)\}$
and
$\{z_{\lambda}^{(\bi)}\:|\:\lambda \in \mathscr P(d), 
\bi \in \Omega(\lambda)\}$ 
give three different bases for
the degree $d$ component of $B$.
Consider the transition matrices
$P = (p_{\lambda,\mu}^{\bi,\bj})$
and 
$Q = (q_{\lambda,\mu}^{\bi,\bj})$
where $\lambda,\mu \in \mathscr P(d)$, $\bi \in \Omega(\lambda), 
\bj \in \Omega(\mu)$
defined from
$$
x_\lambda^{(\bi)} = \sum_{\mu \in \mathscr P(d),\ \bj \in \Omega(\mu)}
p_{\lambda,\mu}^{\bi,\bj} y_\mu^{(\bj)},
\qquad
z_\lambda^{(\bi)} = \sum_{\mu \in \mathscr P(d),\ \bj \in \Omega(\mu)}
q_{\lambda,\mu}^{\bi,\bj} y_\mu^{(\bj)}.
$$

\begin{Lemma}\label{f1} 
The matrix $Q$ is block diagonal, i.e.
$q_{\la,\mu}^{\bi,\bj} = 0$ for $\la \neq \mu$.
Moreover, the determinant of the $\la$-block 
$Q_\la = (q^{\bi,\bj}_{\la,\la})_{\bi,\bj \in \Omega(\la)}$ of $Q$
is $\alpha^{a_\lambda}\beta^{b_\lambda}$,
notation as in the introduction.
\end{Lemma}

\begin{proof}
%
It is immediate from the definition that
$q_{\la,\mu}^{\bi,\bj} = 0$ for $\la \neq \mu$.
So consider the $\la$-block $Q_\la$ of $Q$.
Represent $\la = (\lambda_1 \geq \dots \geq \lambda_h > 0)$
instead as $(1^{r_1} 2^{r_2} \dots s^{r_s})$.
By definition,
$$
z_\la^{(\bi)} = \sum_{\bj \in I(\lambda)}
a^{(\lambda_1)}_{i_1,j_1} \dots a^{(\lambda_h)}_{i_h,j_h} y_{\la}^{(\bj)}.
$$
Thus $Q_\lambda$ is the matrix
$S^{r_1} (A^{(1)}) \otimes \dots \otimes S^{r_s}(A^{(s)})$,
a tensor product of symmetric powers of the matrices $A^{(n)}$.
Now, note that for an $n \times n$ matrix $A$,
$$
\det S^m(A) = (\det A)^{\binom{n+m-1}{n}},
$$
while for an $n \times n$ matrix $B$ and an $m \times m$ matrix $C$,
$$
\det(B \otimes C) = (\det B)^m (\det C)^n.
$$
These are both proved by reducing to the case that the matrices are diagonal.
Combining the formulae with Lemma~\ref{f0}, one computes
$\det Q_\lambda = \alpha^{a_\lambda} \beta^{b_\lambda}$.
\end{proof}

Now we can prove the main theorem:

\vspace{1mm}
\begin{Theorem}
The determinant of the restriction of the Shapovalov form to the
degree $d$ part of $B_\Z$ is
$\prod_{\lambda \in \mathscr P(d)}\alpha^{a_\lambda} \beta^{b_\lambda}$.
\end{Theorem}

\begin{proof}
Consider the matrices
$M = (m_{\la,\mu}^{\bi,\bj})$
and
$N = (n_{\la,\mu}^{\bi,\bj})$
for $\la,\mu \in \mathscr P(d)$, $\bi \in \Omega(\la)$, $\bj \in \Omega(\mu)$
defined from
$$
m_{\la,\mu}^{\bi,\bj} = (x_\la^{(\bi)}, x_\mu^{(\bj)})_S,
\qquad
n_{\la,\mu}^{\bi,\bj} = (x_\la^{(\bi)}, x_\mu^{(\bj)})_K.
$$
Recalling the transition matrices $P$ and $Q$ introduced above,
Lemma~\ref{f2} gives at once that
$M = P Q P^{-1} N.$
On the other hand, $N$ has determinant $1$,
since as we observed above 
the degree $d$ component of 
$B_\Z$ admits an  orthonormal basis with respect to the 
contravariant form $(.\,,.)_K$ (see \cite[Corollary 2.1]{CKK}).
So we can compute $\det M$ at once using Lemma~\ref{f1}.
\end{proof}

We finally indicate how to deduce the generating functions
$a(q)$ and $b(q)$ stated in the introduction.
By definition, $a_\lambda$ is  $$
\frac{h}{\ell}
\times \left(
\begin{array}{l}
\hbox{the number of ways of coloring}\\
\hbox{the parts $\la_i \equiv 0 \:(r)$ with $\ell$ colors}
\end{array}
\right)
\times
\left(
\begin{array}{l}
\hbox{the number of ways of coloring}\\
\hbox{the parts $\la_i \not\equiv 0 \:(r)$ with $k$ colors}
\end{array}
\right)
$$
where $h$ is the number of $\la_i \equiv 0\:(r)$.
Consider
$$
G(q,t,u) = 
\left(\prod_{n \geq 1} \frac{1}{1-q^{nr}t}\right)^{\ell}
\left(\prod_{n \geq 1} \frac{1-q^{nr} u}{1-q^{n}u}
\right)^k.
$$
The coefficient of $q^d t^h u^i$ is equal to the number of 
partitions of $d$ with $h$ parts divisible by $r$ colored by $\ell$
different colors and with $i$ parts not divisible by $r$ colored
by $k$ different colors.
Hence the generating function $a(q)$ 
for $a(d) = \sum_{\lambda \in\mathscr P(d)} a_\lambda$ is equal to 
$\frac{1}{\ell} \frac{d}{dt} G(q,t,u)|_{t=u=1}$.
Similarly, $b(q) = \frac{1}{k} \frac{d}{du} G(q,t,u)|_{t=u=1}$.

\ifx\undefined\bysame
\newcommand{\bysame}{\leavevmode\hbox to3em{\hrulefill}\,}
\fi

\end{document}